\documentclass[11pt]{amsart}
\usepackage{graphicx}
\usepackage{amssymb}
\usepackage{amsmath}
\usepackage{amsthm}
\usepackage{amscd}
\newtheorem{thm}{Theorem}[section]

\newtheorem{cor}[thm]{Corollary}
\newtheorem{lem}[thm]{Lemma}

\theoremstyle{definition}

\theoremstyle{remark}

\numberwithin{equation}{section}

\newcommand{\gr}{{\rm{gr}}}

\begin{document}

\title[On Graded Bialgebra Deformations]{On Graded Bialgebra Deformations}
\author[Y. Du, \ X.W. Chen \ and \   Y. Ye
] {Yu Du, \ Xiao-Wu Chen $^*$  \ and \ Yu Ye}

\thanks{$^*$ The corresponding author}
\thanks{Supported in part by the National Natural Science Foundation of China (Grant No. 10271113 and No. 10501041)
and the Europe Commission AsiaLink project  ``Algebras and
Representations in China and Europe$"$ ASI/B7-301/98/679-11}
\thanks{E-mail: ydu@mail.ustc.edu.cn,
xwchen@mail.ustc.edu.cn,
 yeyu@ustc.edu.cn}

\maketitle

\begin{center}
Department of Mathematics \\University of Science and Technology of
China \\Hefei 230026, Anhui, P. R. China
\vskip 5pt
USTC Shanghai Institute for Advanced Studies\\
Shanghai, 201315, P.R.China
\end{center}

\begin{abstract}
We introduce the graded bialgebra deformations, which explain
Andruskiewitsch-Schneider's liftings method. We also relate this
graded bialgebra deformation with the corresponding graded bialgebra
cohomology groups, which is the graded version of the one due to
Gerstenhaber-Schack.
\end{abstract}

\section{ Introduction}
The classification of finite-dimensional pointed Hopf algebras is a
basic problem in the theory of Hopf algebras. It is well-known that
any pointed Hopf algebra $H$ has a coradical filtration, with
respect to which one associates a coradically-graded Hopf algebra
$\gr H$. Following Andruskiewitsch and Schneider, the classification
problem can be divided into two parts. One is the classification of
all coradically-graded pointed Hopf algebras. The other is to find
all possible pointed Hopf algebras $H$ with $\gr H$ isomorphic to a
given coradically-graded pointed Hopf algebra. The second part is
just the lifting method in \cite{AS} and \cite{AS2}. One of our
motivations is to relate the lifting method with certain bialgebra
deformation theory.

\par \vskip 5pt
The deformation theory for algebras is initiated by Gerstenhaber in
\cite{G}, and its analogue for bialgebras appeared first in
\cite{GS} (also see \cite{GS2} and \cite{PW}). Inspired by the
graded algebra deformation theory in \cite{SS} and \cite{BG}, we
develop in this paper the theory of graded bialgebra deformations
and their corresponding cohomology groups. Moreover this deformation
theory can be used to explain Andruskiewitsch-Schneider's lifting
method.

\par \vskip 5pt
The paper is organized as follows. In section 2, first we recall the
notion of liftings and introduce the graded bialgebra deformations,
and we show that the lifting is just the same as the graded
bialgebra deformation in the sense of Theorem 2.2. The graded-rigid
bialgebras are also studied, see Corollary 2.3 and Corollary 2.4. In
section 3, we introduce the notion of graded ``hat'' bialgebra
cohomology groups for graded bialgebras, which controls the graded
bialgebra deformations, see Theorem 3.3.
\par \vskip 20pt

\section{Liftings and graded bialgebra deformations}
We will work on a base field $\mathbb{K}$. All unadorned tensors are
over $\mathbb{K}$. We refer the notion of  graded bialgebras and
filtered bialgebras to \cite{Sw}, the notion of graded linear maps
to \cite{NV} and \cite{M}.

\subsection{}Let us recall Andruskiewitsch-Schneider's liftings
method, for more details, see \cite{AS2}. Note that the lifting
defined here is a slight generalization.
\par\vskip 5pt

 Throughout, $B=\oplus_{i \geq 0} B_{(i)}$ will be a  graded
bialgebra over $\mathbb{K}$, with identity element $1_B$,
multiplication map $m$, counit $\varepsilon$, and comultiplication
$\Delta$. Then $B$ has a natural bialgebra filtration
\begin{align*}
B_0 \subseteq B_1 \subseteq \cdots \subseteq B_i \subseteq \cdots,
\end{align*}
where $B_i=\oplus_{j \leq i}B_{(j)}$ for any $n \geq 0$.\par \vskip
5pt

A \emph{lifting} of the graded bialgebra $B$ is a filtered bialgebra
structure, denoted by $U$,  on the underlying filtered vector space
$B$ with the above filtration such that
\begin{align*}
{\rm gr}U =B
\end{align*}
as graded bialgebras, where ${\rm gr}U$ is the graded bialgebra
associated to the filtered bialgebra $U$ (\cite{Sw}, p.226). (By
${\rm gr}U =B$,  we use the natural identification of the underlying
space ${\gr U}$ with $B$, that is ${\rm gr} U_{(i)}=B_i/B_{i-1}\cong
B_{(i)}$ for each $i \geq 0$.)\par \vskip 5pt

For any lifting $U$ of the graded bialgebra $B$, it follows from the
definition that $U$ and $B$ have the same identity element and the
counit. Therefore,  to give a lifting $U$, we just need to define
the multiplication $m_U$ and comultiplication $\Delta_U$ .\par

Two liftings $U$, $V$ of the graded bialgebra $B$ are said to be
\emph{equivalent}, if there is filtered bialgebra isomorphism
$\theta: U \longrightarrow V$ such that ${\rm gr \theta}={\rm
Id}_B$, where ${\rm gr} \theta$ is the graded morphism associated
to $\theta$, and here again we use the identifications ${\gr U}=B$
and ${\gr V}=B$ (as graded bialgebras).
\par \vskip 5pt

 Denote by
\begin{align*}
\emph{Lift}(B)
\end{align*}
 the set of equivalent classes of all the liftings of the graded
bialgebra $B$.

\subsection{}
In this subsection, we will study graded bialgebra deformations of
the graded bialgebra $B=\oplus _{i \geq 0} B_{(i)}$.\par \vskip 5pt

Let $l\in \mathbb{N} \cup \{+ \infty\}$. Consider the space
$B[t]/{(t^{l+1})}$, which is viewed as a free module over
$\mathbb{K}[t]/{(t^{l+1})}$, and also a graded $\mathbb{K}$-space
with ${\rm deg }t=1$ and ${\rm deg} b=n$, if $b \in B_{(n)}$. If
$l=+\infty$, then $B[t]/{(t^{l+1})}$ means $B[t]$ and
$\mathbb{K}[t]/{(t^{l+1})}$ means $\mathbb{K}[t]$.
\par \vskip 10pt

An  \emph{$l$-th level graded bialgebra deformation} of $B$ consists
of
\begin{align*}
m_t^l: (B\otimes B)[t]/{(t^{l+1})} \longrightarrow B[t]/{(t^{l+1})}
\end{align*}
and \begin{align*}
 \Delta_t^l: B[t]/{(t^{l+1})} \longrightarrow
(B\otimes B)[t]/{(t^{l+1})}\cong
B[t]/{(t^{l+1})}\otimes_{\mathbb{K}[t]/{(t^{l+1})}} B[t]/{(t^{l+1})}
\end{align*}
 which are $\mathbb{K}[t]/{(t^{l+1})}$-linear and homogeneous maps of
degree zero such that
\par \vskip 3pt

\begin{enumerate}
\item[(\textit{i})] $B[t]/{(t^{l+1})}$ is a bialgebra over
$\mathbb{K}[t]/{(t^{l+1})}$ with identity element $1_B$,
multiplication $m_t^l$, counit $\varepsilon_t^l$ and
comultiplication $\Delta_t^l$, where the counit $\varepsilon_t^l:
B[t]/{(t^{l+1})} \longrightarrow \mathbb{K}[t]/{(t^{l+1})}$ is
given by $\varepsilon_t^l(bt^j)= \varepsilon(b)t^j$, $b \in B$, $0
\leq j \leq l$;\par \vskip3 pt

\item[(\textit{ii})] $m_t^l \equiv m\otimes {\rm
Id}_{\mathbb{K}[t]/{(t^{l+1})}}$ and  $\Delta_t^l \equiv
\Delta\otimes {\rm Id}_{\mathbb{K}[t]/{(t^{l+1})}}$ \ mod$(t)$,
where $m$ and $\Delta$ are the multiplication and comultiplication
of $B$, respectively.
\end{enumerate}
\par \vskip 5pt

 Denote by $(B[t]/{(t^{l+1})}, m_t^l,
\Delta_t^l)$ the above $l$-th level graded bialgebra deformation.
\par \vskip 10pt

From now on, we will abbreviate $l$-th level graded bialgebra
deformations as $l$-deformations, and ${+ \infty}$-deformations will be referred simply as deformations.
 Denote by $\mathcal{E}^l(B)$ the
set of all $l$-deformations of the graded bialgebra $B$,
 and $\mathcal{E}^{+ \infty}(B)$ is written as $\mathcal{E}(B)$.
 Elements of $\mathcal{E} (B)$ will be written as $(B[t], m_t,
 \Delta_t)$.
 \par \vskip
5pt

Two $l$-deformations $(B[t]/{(t^{l+1})}, m_t^l, \Delta_t^l)$ and
$(B[t]/{(t^{l+1})}, {m'}_t^l, {\Delta'}_t^l) $ are said to be
\emph{isomorphic}, if there exists an isomorphism of
$\mathbb{K}[t]/{(t^{l+1})}$-bialgebras
\begin{align*}
\phi: (B[t]/{(t^{l+1})}, m_t^l, \Delta_t^l) \longrightarrow
(B[t]/{(t^{l+1})}, {m'}_t^l, {\Delta'}_t^l)
\end{align*}
such that   $\phi$ is homogeneous of degree zero and
\begin{align*}\phi \equiv {\rm Id}_B \otimes {\rm
Id}_{\mathbb{K}[t]/{(t^{l+1})}} \ \mbox{mod}(t).
\end{align*}
\par \vskip 5pt

Denote by
\begin{center}
$iso \mathcal{E}^l(B)$ \ (\emph{resp}. \ $iso \mathcal{E}(B)$)
\end{center}
the set of isoclasses of $l$-deformations (\emph{resp}.
deformations) of the graded bialgebra $B$, for $l \in \mathbb{N}$.
\par \vskip 5pt

\subsection{}
Use the notation as above. Consider an element $(B[t]/{(t^{l+1})},
m_t^l, \Delta_t^l)$ of $\mathcal{E}^l(B)$. By definition, we can
write
\begin{align}
m_t^l(a \otimes b)&=\sum_{0 \leq s \leq l} m_s(a \otimes b) t^s,
\end{align}
and
\begin{align}
\Delta_t^l(c)&=\sum_{0 \leq s \leq l} \Delta_s(c)t^s,
\end{align}
where $a, b, c\in B$, and $m_s:B \otimes B \longrightarrow B$ and
$\Delta_s:B \longrightarrow B\otimes B$ are homogeneous of degree
$-s$. Note that $m_0=m$ and $\Delta_0=\Delta$. \par \vskip 5pt

It is easy to check that the associativity of $m_t^l$, the
compatibility of $m_t^l$ and $\Delta_t^l$, and the coassociativity
of $\Delta_t^l$ are equivalent to the following identities,
respectively, for each $1\leq n \leq l$,

\begin{align}
&a  m_n(b\otimes c)-m_n(ab \otimes c)+ m_n(a\otimes
bc)-m_n(a\otimes
b)c\\
=&\sum_{1 \leq s \leq n-1} m_s(m_{n-s}(a\otimes b)\otimes
c)-m_s(a\otimes m_{n-s}(b\otimes c)), \nonumber
\end{align}

\begin{align}
&m_n(a_{(1)}\otimes b_{(1)})\otimes
a_{(2)}b_{(2)}-\Delta(m_n(a\otimes b))+ a_{(1)}b_{(1)}\otimes
m_n(a_{(2)}\otimes b_{(2)})\\
 &\quad + a_{(1)}b_l \otimes a_{(2)}b_r -\Delta_n (ab) +a_lb_{(1)}\otimes a_rb_{(2)} \nonumber \\
=&-\sum_{0\leq s,r,s',r'\leq n-1, \;
 s+s'+r+r'=n} (m_r \otimes
m_{r'})\circ \tau_{23}\circ (\Delta_s\otimes \Delta_{s'})(a\otimes
b)\nonumber \\
&+ \sum_{1 \leq s\leq n-1} \Delta_s(m_{n-s}(a \otimes
b)) ,\nonumber
\end{align}
and

\begin{align}
&c_{(1)}\otimes \Delta_n(c_{(2)})- (\Delta\otimes {\rm Id})\circ
\Delta_n(c)  +({\rm Id}\otimes \Delta)\circ \Delta_n(c)
-\Delta_n(c_{(1)})\otimes
c_{(2)}\\
=&\sum_{1 \leq s \leq n-1} (\Delta_{n-s} \otimes {\rm Id})\circ
\Delta_s(c)-({\rm Id}\otimes \Delta_{n-s})\circ \Delta_s(c)
,\nonumber
\end{align}
where we use Sweedler's notation $\Delta(a)=a_{(1)}\otimes
a_{(2)}$, $a \in B$, and  in the second identity we use the
notation $\Delta_n(a)=a_l\otimes a_r$ and $\Delta_n(b)=b_l\otimes
b_r$, and the map $\tau_{23}$ is the canonical flip map at the
second and third positions.
\par \vskip 5pt

Let $(B[t]/{(t^{l+1})}, m_t^l, \Delta_t^l)$ and $(B[t]/{(t^{l+1})},
{m'}_t^l, {\Delta'}_t^l)$ be two $l$-deformations with the maps
$m_s$, $\Delta_s$ and $m'_s$, $\Delta'_s$ as in (2.1) and (2.2). An
isomorphism $\phi$ between these deformations is given by
\begin{align}
\phi(a)=\sum_{0 \leq s \leq l} \phi_s(a)t^s, \quad a \in B,
\end{align}
where $\phi_s:B \longrightarrow B$ is a homogeneous map of degree
$-s$. Note that $\phi_0={\rm Id}_B$. The fact that $\phi$ is a
morphism of $\mathbb{K} [t]/{(t^{l+1})}$-bialgebras implies $\phi$
preserves the identity element $1_B$ and the counit
$\varepsilon_t^l$, and it satisfies, for each $ 1 \leq n \leq l$,
\begin{align}
&(m_n-m'_n)(a\otimes b)= a\phi_n(b)-\phi_n(ab)+\phi_n(a)b\\
           &\quad +\sum_{0<s<n} \{ \phi_s(a)\phi_{n-s}(b) -\phi_{s}(m_{n-s}(a\otimes
           b))+\sum_{r+r'=n-s} m'_s(\phi_r(a)\otimes
           \phi_{r'}(b))\} \nonumber
\end{align}

and

\begin{align}
&(\Delta_n-\Delta'_n)(c)= \Delta(\phi_n(c))-c_{(1)}\otimes
\phi_n(c_{(2)})
-\phi_n(c_{(1)})\otimes c_{(2)} \\
&\quad + \sum_{0< s< n} \{ \Delta'_s(\phi_{n-s}(c)) -(\phi_s \otimes
\phi_{n-s})(\Delta(c))
 - \sum_{r+r'=n-s} (\phi_r\otimes \phi_{r'})
(\Delta_s(c))\} \nonumber,
\end{align}
for all $a, b, c \in B$.
 Note that above discussion works for
all $l \in \mathbb{N} \cup \{+\infty\}$.\par \vskip 5pt

The analogue of the following lemma is well-known in classical
deformation theory.

\begin{lem}There exist  restriction maps $r_{l, l'}: \mathcal{E}^l(B) \longrightarrow \mathcal{E}^{l'}(B)$
for every $l> l' \in \mathbb{N}$, and maps $r_{l}: \mathcal{E}(B)
\longrightarrow \mathcal{E}^l(B)$ such that
\begin{align*}
 \mathcal{E}(B) = \underleftarrow{\rm lim}_{l \in \mathbb{N}} \   \mathcal{E}^l(B).
\end{align*}
\end{lem}

\noindent{\bf Proof.} The restriction map $r_{l,l'}$ is given as
follows: given $(B[t]/{(t^{l+1})}, m_t^l, \Delta_t^l)$ in
$\mathcal{E}^l(B)$ with the maps $m_s$ and $\Delta_s$ defined in
(2.1) and (2.2), just define $m_t^{l'}:=\sum_{0 \leq s \leq l'} m_s
t^s$ and $\Delta_t^{l'} :=\sum_{0 \leq  s \leq l'} \Delta_s t^s$; it
is direct to check that $(B[t]/{(t^{l'+1})}, m_t^{l'},
\Delta_t^{l'})$ is the desired element in $\mathcal{E}^{l'}(B)$. The
map $r_l$ is defined in a similar way, and then the result is
obvious. \hfill $\blacksquare$

           \par \vskip 10pt

   A graded bialgebra $B=\oplus_{i\geq 0} B_{(i)}$ is called \emph{graded-rigid} if the set $iso \mathcal{E}(B)$ has
   only one element, i.e., any deformation of $B$ is isomorphic to the trivial
   one.

\subsection{}
We have the following observation, which says that the graded
bialgebra deformations coincide with the liftings.

\begin{thm}
Let $B=\oplus_{i \geq 0} B_{(i)}$ be a graded bialgebra. There
exists a natural bijection
\begin{align*}
Lift(B) \simeq \ iso \mathcal{E}(B).
\end{align*}
\end{thm}

\noindent{\bf Proof.}\quad  We will construct a map $F: Lift(B)\longrightarrow \ iso \mathcal{E}(B)$. Given a lifting $U$
of $B$.
Denote by $m_U$ and $\Delta_U$ the multiplication and comultiplication maps of $U$.
 Since $U$ is a filtered bialgebra, we have
\[
m_U: B_i\otimes B_j \longrightarrow B_{i+j} \quad \text{and} \quad
 \Delta_U: B_n \longrightarrow \sum_{i+j=n} B_i\otimes B_j.
\]
Therefore, for any $s\ge0$, there uniquely exist homogeneous maps of
degree $-s$, say $m_s: B\otimes B \longrightarrow B$ and $\Delta_s:B
\longrightarrow B\otimes B$, such that
\[
m_U(a \otimes b)=\sum_{s \geq 0} m_s(a \otimes b) \quad \text{and}
 \quad \Delta_U(c)=\sum_{s \geq 0} \Delta_s(c).
\]
By ${\rm gr}U=B$ as graded bialgebras, we have $m_0=m$ and
$\Delta_0=\Delta$. \par \vskip 5pt

 Now Define $F(U)=(B[t], m_t,
\Delta_t)$ as follows
  \begin{align*}
m_t(a \otimes b):=\sum_{s \geq 0} m_s(a \otimes b) t^s \quad
\mbox{and }\quad  \Delta_t(c):=\sum_{s \geq 0} \Delta_s(c)t^s.
\end{align*}
It is direct to check that $F(U)$ is a deformation.

$F$ is well-defined, i.e., it maps equivalent liftings to isomorphic
deformations. In fact, for given liftings $U$ and $V$, an
equivalence $\theta$ of $U$ and $V$ is a filtered isomorphism, hence
for any $s\ge0$, there determines a unique homogeneous map $\phi_s:
B \longrightarrow B$ of degree $-s$ such that
\begin{align*}
\theta(a)= \sum_{s\geq 0} \phi_s(a), \quad  a \in B.
\end{align*}
Then define a $\mathbb{K}[t]$-linear map $\phi:B[t]\longrightarrow
B[t]$ such that $\phi(a)=\sum_{s \geq 0}\phi_s(a)t^s $. Hence $\phi$
is an isomorphism between the deformations $F(U)$ and $F(V)$.\par

 On the other hand, by seeing (2.1) and (2.2), one
obtains that $F$ is a bijection. This completes the proof. \hfill
$\blacksquare$

                \vskip 20pt

An immediate consequence of  Theorem 2.2 is
\begin{cor}
Let $B=\oplus_{i \geq 0}B_{(i)}$ be a graded bialgebra. Then  $B$ is
graded-rigid implies that, for any filtered bialgebra $U$ such that
${\rm gr}U\simeq B$ as graded bialgebras, we have $U\simeq B$ as
bialgebras.\par

 If we assume that $B$ is coradically-graded, the
converse is also true.
\end{cor}

\noindent{\bf Proof.}\quad By Theorem 2.2, $B$ is graded-rigid if
and only if $Lift(B)$ is a single element set, i.e., every lifting
of $B$ is trivial. \par \vskip 5pt

For the first statement, such a filtered bialgebra $U$ with ${\rm
gr}U\simeq B$ gives rise to a lifting on $B$, denoted by $U'$, such
that $U \simeq U'$ (as bialgebras). Since $B$ is graded-rigid, we
get $U'\simeq B$, thus we are done. \par \vskip 3pt

For the second one, assume $B$ is coradically-graded. Let $U$ be a
lifting of $B$. Thus by the assumption, there exists an isomorphism
$\theta:U \simeq B$. Note that $\theta$ preserves the coradical
filtration, thus ${\rm gr} \theta$ can be viewed as a graded
automorphism of $B$. Thus take $\theta'= ({\rm gr}\theta)^{-1}\circ
\theta: U \simeq B$. So $\theta'$ realizes an equivalence between
the lifting $U$ and the trivial lifting. This proves that $B$ is
graded-rigid. \hfill $\blacksquare$

\subsection{}
In this subsection, we assume that the base field $\mathbb{K}$ is
algebraically closed of characteristic zero. One can define the
variety ${\rm Bialg}_n$ of the bialgebra structures on
$n$-dimensional spaces, which carries a natural
$GL_n(\mathbb{K})$-action by base changes, see \cite{St} and
\cite{Ma}. Recall that a bialgebra $B$ is called rigid if
$GL_n(\mathbb{K})$-orbit of ${\rm Bialg}_n$ containing $B$ is
Zariski open. In fact, we have

\begin{cor}
Let $\mathbb{K}$ be an algebraically closed field of characteristic
zero, $B=\oplus_{i \geq 0}B_{(i)}$ a finite dimensional graded
bialgebra over $\mathbb{K}$. If  $B$ is rigid and
coradically-graded, then $B$ is graded-rigid in the sense of {\rm
\bf 2.3}.
\end{cor}

\noindent{\bf Proof.}\quad By Corollary 2.3, we only need to show
that every filtered bialgebra $U$ with ${\rm gr}U\simeq B$ is
isomorphic to $B$. Assume the dimension of $B$ is $n$. By Theorem
3.4 in \cite{Ma}, $B$ is a degeneration of $U$, i.e., lies the
closure of the orbit of $U$ (in the variety ${\rm Bialg}_n$).
However the $GL_n(\mathbb{K})$-orbit of $B$ is open, we obtain that
$B$ and $U$ belong to the same $GL_n(\mathbb{K})$-orbit, i.e.,  $B
\simeq U$ as bialgebras, finishing the proof.\hfill $\blacksquare$

\par \vskip 20pt

\section{Graded bialgebra cohomology}
 In this section we will relate  the graded bialgebra deformations with corresponding
 cohomology groups, which
 will be a graded (and normalized) version of ``hat"  bialgebra cohomology groups introduced in
 \cite{GS}
 (also see \cite{PW}).

\subsection{} Let $(B, m, e, \Delta, \varepsilon)$ be a bialgebra.
Again we will use Sweedler's notation $\Delta(a)=a_{(1)}\otimes
a_{(2)}$, $a\in B$.\par \vskip 5pt
 Let us recall the bicomplex in \cite{GS} or
\cite{PW}, p.619. For this end, we need the following maps, where
$p, q \geq 1$ and all $b$'s are in $B$, $\lambda^p: B^{\otimes
{p+1}} \longrightarrow B^{\otimes p}$ and $\rho^p: B^{\otimes {p+1}}
\longrightarrow B^{\otimes p}$ are given by
\begin{align*}
\lambda^p(b^1 \otimes \cdots \otimes b^{p+1})&= b^{1}_{(1)}b^2 \otimes \cdots \otimes
b^1_{(p)}b^{p+1},\\
\rho^p(b^1 \otimes \cdots \otimes b^{p+1})&= b^1 b^{p+1}_{(1)}\otimes \cdots  \otimes b^p b^{p+1}_{(p)}.
\end{align*}
Dually, the maps $\sigma^q: B^{\otimes q} \longrightarrow B^{\otimes {q+1}}$ and
$\tau^q: B^{\otimes q} \longrightarrow B^{\otimes q+1}$ are given by
\begin{align*}
\sigma^q(b^1 \otimes \cdots \otimes b^q)&=(b^1_{(1)} \cdots b^q_{(1)}) \otimes b^1_{(2)} \otimes \cdots b^q_{(2)},
\\
\tau^q(b^1 \otimes \cdots \otimes b^q)&= b^1_{(1)} \otimes \cdots \otimes b^q_{(1)} \otimes (b^1_{(2)} \cdots b^q_{(2)}).
\end{align*}
In addition, we need $\Delta_i^p: B^{\otimes p} \longrightarrow B^{\otimes {p+1}}$ and $\mu_j^q: B^{\otimes q+1}
\longrightarrow B^{\otimes q}$, $1 \leq i \leq p$ and $1 \leq j \leq q $, which are given by
\begin{align*}
&\Delta_i^p(b^1 \otimes \cdots \otimes b^p)=b^1 \otimes \cdots \otimes b^i_{(1)} \otimes b^i_{(2)} \otimes \cdots \otimes
b^p,\\
&\mu_i^q(b^1 \otimes \cdots \otimes b^{q+1}) =b^1 \otimes \cdots \otimes b^ib^{i+1} \otimes \cdots \otimes b^{q+1}.
\end{align*}

\par \vskip 8pt

Let $C^{p, q}={\rm Hom}_\mathbb{K}(B^{\otimes q}, B^{\otimes p})$,
$p, q \geq 1$. Define
\begin{align*}
 \delta_h^{p, q}: C^{p, q} \longrightarrow C^{p, q+1} \quad \mbox{ and }
\quad  \delta_c^{p, q}: C^{p, q} \longrightarrow C^{p+1, q}
\end{align*}
which are given by
\begin{align*}
\delta_h^{p, q} (f)&= \lambda^p \circ ({\rm Id}\otimes f) + \sum_{i=1}^q (-1)^i f \circ \mu_i^q + (-1)^{q+1} \rho^p \circ (f \otimes {\rm
Id})\\
\delta_c^{p, q} (f)&= ({\rm Id} \otimes f) \circ \sigma^q +\sum_{j=1}^p (-1)^j \Delta_j^p \circ f + (-1)^{p+1} (f \otimes {\rm Id}) \circ \tau^q
\end{align*}
 for $f \in C^{p, q}$, where ${\rm Id}$ denotes the identity map of
 $B$.\par \vskip 5pt

 It is direct to check that $(C^{p, q}, \delta_h^{p, q}, \delta_c^{p,
 q})$ is  a bicomplex (see \cite{PW}, p.619), i.e.,
 \begin{align*}
 \delta_h^{p, q+1} \circ \delta_h^{p, q}=0, \quad  \delta_c^{p, q+1} \circ \delta_h^{p,
 q}=\delta_h^{p+1, q} \circ \delta_c^{p,q}, \quad \delta_c^{p+1, q}\circ
 \delta_c^{p, q}=0.
 \end{align*}

 \par \vskip 5pt
We will introduce a sub-bicomplex of the above bicomplex. Let
$\textsf{m}={\rm Ker}\varepsilon$.  Denote by $\textsl{i}:
\textsf{m} \longrightarrow B$ the inclusion map, and $\pi: B
\longrightarrow \textsf{m}$ is given by  $\pi(b)= b
-\varepsilon(b)1_B$,  $b \in B$.  Set $D^{p, q}= {\rm
Hom}_\mathbb{K}(\textsf{m}^{\otimes q}, \textsf{m}^{\otimes p})$,
$p, q \geq 1$. Note that we have a natural embedding  $D^{p, q}
\hookrightarrow C^{p, q}$ by identifying
\begin{align*}
f \in D^{p, q}  \quad \mbox{   with   }  \quad \textsl{i}^{\otimes p} \circ f \circ \pi^{\otimes q} \in C^{p,q}.
\end{align*}

\par \vskip 5pt
We have the following observation

\begin{lem}
Use the above notation.
 Then $\delta_h^{p, q}(D^{p, q}) \subseteq D^{p, q+1}$ and $\delta_c^{p, q} (D^{p, q}) \subseteq D^{p+1, q}$.
\end{lem}

\noindent {\bf Proof.}\quad  Just note that $f \in C^{p, q}$ lies in
$D^{p, q}$ if and only if
\begin{center}
$({\rm Id}^{\otimes j-1} \otimes \varepsilon \otimes {\rm
Id}^{\otimes p-j}) \circ f=0$
\end{center}
and
 \begin{center}

 $f(b^1\otimes \cdots \otimes b^{i-1} \otimes  1\otimes b^{i+1}
\otimes \cdots \otimes b^q )=0$, \end{center}

  for any $1 \leq i
\leq q$, $1 \leq j \leq p$ and any $b^i\in B$.   Then the lemma
follows from the definition of $\delta_h^{p, q}$ and $\delta_c^{p,
q}$ immediately. \hfill $\blacksquare$

\subsection{}
From now on $B=\oplus_{i \geq 0} B_{(i)}$ will be  a graded
bialgebra. In this case $\textsf{m} \subseteq B$ is a graded
subspace. Consider $D^{p, q}_{(l)}:= {\rm
Hom}_\mathbb{K}(\textsf{m}^{\otimes q}, \textsf{m}^{\otimes
p})_{(l)}$ , $l \in \mathbb{Z}$, whose elements are homogeneous maps
from $\textsf{m}^{\otimes q}$ to $\textsf{m}^{\otimes p}$ of degree
$l$. Note that $D^{p, q}_{(l)} \subseteq D^{p, q} \hookrightarrow
C^{p, q}$. We have the following

\begin{lem}
$\delta_h^{p, q}(D^{p, q}_{(l)}) \subseteq D^{p, q+1}_{(l)}$ and
$\delta_c^{p, q}(D^{p, q}_{(l)}) \subseteq D^{p+1, q}_{(l)}$ for
each $l \in \mathbb{Z}$, $p, q \geq 1$.
\end{lem}

\noindent {\bf Proof.}\quad Set $C^{p, q}_{(l)}= {\rm
Hom}_\mathbb{K}(B^{\otimes q}, B^{\otimes p})_{(l)}$.  Clearly,
$D^{p, q}_{(l)}= D^{p, q} \cap C^{p, q}_{(l)}$. From the definition
of $\delta_h^{p, q}$  and $\delta_c^{p, q}$, one sees that they
preserve the degrees, i.e.,
 $\delta_h^{p, q}(C^{p, q}_{(l)})
\subseteq C^{p, q+1}_{(l)}$ and  $\delta_c^{p, q}(C^{p, q}_{(l)})
\subseteq C^{p+1, q}_{(l)}$. Now the result follows from  Lemma 3.1.
\hfill $\blacksquare$
\par \vskip 10pt

Denote by $\delta_{h, (l)}^{p, q}$ (\emph{resp}. $\delta_{c,
(l)}^{p, q}$) the restriction of the maps $\delta_h^{p, q}$
(\emph{resp}. $\delta_c^{p, q}$) to the subspace $D^{p, q}_{(l)}$.
Thus by Lemma 3.2, we get a bicomplex $(D_{(l)}^{p, q}, \delta_{h,
(l)}^{p, q}, \delta_{c, (l)}^{p, q})$ for each $l \in
\mathbb{Z}$.\par
 There is a canonical way to construct a complex
from a given bicomplex. Set
 \begin{align*}
 \hat{D}^n_{(l)}= \bigoplus_{p+q=n+1, p,
q\geq 1} D_{(l)}^{p, q}, \quad n \geq 1;
\end{align*}
define $\partial_{(l)}^n : \hat{D}^n_{(l)} \longrightarrow
\hat{D}^{n+1}_{(l)} $
 by
\begin{align*}
\partial^n_{(l)}|_{D^{n+1-q, q}_{(l)}}:= \delta_{h, (l)}^{p, q} +
(-1)^q \delta_{c, (l)}^{p, q}, \quad 1 \leq q \leq n.
\end{align*}
Hence, for each $l \in \mathbb{Z}$, we get a complex

\begin{align*}
0 \longrightarrow \hat{D}^1_{(l)}
\stackrel{\partial_{(l)}^1}{\longrightarrow} \hat{D}^2_{(l)}
\stackrel{\partial_{(l)}^2}{\longrightarrow} \hat{D}^3_{(l)}
\stackrel{\partial_{(l)}^3}{\longrightarrow} \hat{D}^4_{(l)}
\longrightarrow \cdots
\end{align*}
\par \vskip 5pt

We define the n-th cohomology group of the above complex to be the
\emph{n-th graded ``hat'' bialgebra cohomology  of degree $l$} of
the graded bialgebra $B$, which will be denoted by
$\hat{h}_b^n(B)_{(l)}$, $n \geq 1$, $l \in \mathbb{Z}$.
\par \vskip 8pt

It is very useful to write out $\hat{h}_b^2(B)_{(l)}$ and
$\hat{h}_b^3(B)_{(l)}$ explicitly from the definition. In what
follows, we will use the maps $\delta_{h}^{p, q}$ and
$\delta_{c}^{p, q}$, instead of $\delta_{h, (l)}^{p, q}$ and
$\delta_{c, (l)}^{p, q}$ for simplicity. We have the following
facts.

\vskip5pt {\bf 1.}\quad The cohomology group $\hat{h}_b^2(B)_{(l)}$
consists of all pairs $(f, g)$, where $f: \textsf{m} \otimes
\textsf{m} \longrightarrow \textsf{m} $ and $g: \textsf{m}
\longrightarrow \textsf{m} \otimes \textsf{m}$ are homogeneous maps
of degree $l$, satisfying the following relations:
\begin{align*}
\delta_h^{1, 2}(f)=0, \quad \delta_c^{1,2}(f)+\delta_h^{2,1}(g)=0,
\quad \delta_c^{2,1}(g)=0,
\end{align*}
i.e., for any $a, b, c \in \textsf{m}$, we have
\begin{align}
&af(b \otimes c)- f(ab \otimes c) + f(a \otimes  bc) -f(a \otimes
b)c=0,\\
&f(a_{(1)} \otimes b_{(1)}) \otimes a_{(2)}b_{(2)}- \Delta(f(a
\otimes b)) + a_{(1)}b_{(1)}\otimes f(a_{(2)} \otimes b_{(2)})  \\
&\qquad + a_{(1)}g(b)_l \otimes a_{(2)} g(b)_r -g(ab)+ g(a)_l
b_{(1)}
\otimes g(a)_r b_{(2)}=0,\nonumber \\
&c_{(1)} \otimes g(c_{(2)})-(\Delta \otimes {\rm Id}) (g(c))+ ({\rm
Id}\otimes \Delta)(g(c))- g(c_{(1)}) \otimes c_{(2)}=0,
\end{align}

where we write  $g(b)=g(b)_l \otimes g(b)_r$,  $b \in B$.
\par
\vskip 8pt

Two pairs $(f, g)=(f', g')$ in $\hat{h}_b^2(B)_{(l)}$ if and only if
there exists a homogeneous map $\theta: \textsf{m} \longrightarrow
\textsf{m}$ of degree $l$ such that, for any $a, b , c \in
\textsf{m}$,
\begin{align}
&(f-f')(a \otimes b)= a \theta(b)-\theta(ab)+ \theta(a)b,\\
&(g-g')(c)= \Delta(\theta(c))-c_{(1)}\otimes \theta(c_{(2)})-
\theta(c_{(1)})\otimes c_{(2)}.
\end{align}

\vskip 5pt {\bf2.}\quad The group $\hat{h}_b^3(B)_{(l)}$ consists of
all triples $(F, H, G)$, where
\begin{align*}
F:\textsf{m}\otimes\textsf{m}\otimes\textsf{m}\longrightarrow
                         \textsf{m}, \quad H: \textsf{m}\otimes
                         \textsf{m}\longrightarrow \textsf{m}\otimes
                         \textsf{m}, \quad G:
                         \textsf{m}\longrightarrow \textsf{m}
                         \otimes \textsf{m}\otimes \textsf{m}
\end{align*}
are homogeneous maps of degree $l$, subject to the relations:
\begin{align*}
\delta_h^{1, 3}(F)=0, \quad \delta_h^{2,2}(F)=\delta_c^{1,3}(H),
\quad \delta_c^{2,2}(H)=-\delta_h^{1,3}(G), \quad
\delta_c^{3,1}(G)=0.
\end{align*}
Note that $(F, H, G)=0$ in $\hat{h}_b^3(B)_{(l)}$ if and only if
there exists $(f, g )\in \hat{D}_{(l)}^{2}$ such that
\begin{align}
(F, H, G)=\partial_{(l)}^2( (f,  g) ),
\end{align}
which can be written out explicitly by the definition of
$\partial_{(l)}^2$.\vskip 10pt

\subsection{} Now we are at the position to present our main
observations, which relate the graded bialgebra deformations of the
graded bialgebra $B$ with the cohomology groups
$\hat{h}_b^2(B)_{(l)}$ and $\hat{h}_b^3(B)_{(l)}$(compare \cite{GS},
Section 5).
\par \vskip 10pt

\begin{thm}Let $B=\oplus_{i \geq 0} B_{(i)}$ be a graded bialgebra.
Use the notation as above. Then \par

\begin{enumerate}
\item[(1).] There is a bijection between $iso \mathcal{E}^1(B)$ and
$ \hat{h}_b^2(B)_{(-1)}$.\par

\item[(2).] If $ \hat{h}_b^2(B)_{(-l)}=0$ for each $l\geq 1$,
then the graded bialgebra $B$ is graded-rigid.\par

\item[(3).] The obstruction to extend an element of
$\mathcal{E}^l(B)$ to $\mathcal{E}^{l+1}(B)$ lies in
$\hat{h}_b^3(B)_{(-l-1)}$, $l \geq 1$. In particular, if
$\hat{h}_b^3(B)_{(-l-1)}=0$, one can extend any element of
$\mathcal{E}^l(B)$ to $\mathcal{E}^{l+1}(B)$.
\end{enumerate}
\end{thm}
\par \vskip 10pt

\noindent {\bf Proof.}\quad (1).\ Recall from  {\bf 2.2} that an
element in $\mathcal{E}^1(B)$ is just given by $(B[t]/{(t^2)},
m_t^1, \Delta_t^1)$. As in {\bf 2.3}, write
\begin{align*}
m_t^1(a \otimes b)= ab + f(a \otimes b)t, \quad
\Delta_t^1(c)=\Delta(c)+ g(c)t,
\end{align*}
where $f: B \otimes B \longrightarrow B$ and $g: B \longrightarrow
B\otimes B$ are homogeneous of degree $-1$. Note that $1_B$ is the
identity element of $B[t]/{(t^2)}$, hence $f(1_B \otimes b)=f(b
\otimes 1_B)=0$ for all $b \in B$. Moreover, for $a, b \in
\textsf{m}$, $\varepsilon_t^1(m_t^1(a\otimes b))=0$ implies that
$\varepsilon_t^1(ab + f(a \otimes b)t)=0$, i.e., $f(a\otimes b) \in
\textsf{m}$. Thus we may view $f$ belongs to $D^{1,2}_{(-1)}$.
Dually one can show that $g \in D^{2,1}_{(-1)}$.

\par \vskip 3pt

Note that $m_t^1$ is an associative multiplication on
$B[t]/{(t^2)}$, thus we get
\begin{align*}
f(a \otimes b)c- f(a\otimes bc) +f(ab\otimes c)- af(b \otimes
c)=0,\quad \forall a, b, c \in B.
\end{align*}
Therefore we get equation (3.1). Similarly, the fact that
$\Delta_t^1$ is an algebra morphism (\textmd{resp}. that
$\Delta_t^1$ is an coassociative comultiplication) gives us equation
(3.2)(\emph{resp}. equation (3.3)), i.e., $(f, g)$ can be viewed as
an element in $\hat{h}_b^2(B)_{(-1)}$.\par \vskip 3pt

Suppose that $(B[t]/{(t^2)}, m_t^1, \Delta_t^1)$ and $(B[t]/{(t^2)},
{m'}_t^1, {\Delta'}_t^1)$ are two isomorphic deformations, with $(f,
g)$ and $(f', g')$ defined as above, respectively. Let $\phi$ (see
also {\bf 2.3}) be the isomorhism. We may write
\begin{align*}
\phi(a)=a + \theta(a)t, \quad \forall a \in B,
\end{align*}
for some homogeneous map $\theta: B \longrightarrow B$ of degree
$-1$ (note that the map  $\theta$ may be viewed as a map from
$\textsf{m}$ to $\textsf{m}$).  Now it is direct to check that
$\theta$ realizes an equivalence of $(f, g)$ and $(f', g')$ in
$\hat{h}_b^2(B)_{(-1)}$. Now we have obtained a map from
$\mathcal{E}^1(B)$ to $\hat{h}_b^2(B)_{(-1)}$, sending
$(B[t]/{(t^2)}, m_t^1, \Delta_t^1)$ to $(f, g)$. One can easily see
that  the correspondence is bijective, as required. \par \vskip 10pt

\noindent (2). \  To prove that $B$ is graded-rigid, we just need to
show that $iso \mathcal{E}(B)$ is a single-element set.\par

Let $(B[t], m_t, \Delta_t)$ be an element in $\mathcal{E}(B)$. As
before, write
\begin{align*}
m_t(a \otimes b)=\sum_{s= 0}^{\infty} m_s(a\otimes b) t^s\quad
\mbox{and} \quad \Delta_t(c)=\sum_{s= 0}^{\infty } \Delta_s(c)t^s.
\end{align*}
Note that $m_0=m$ and $\Delta_0=\Delta$, and  $m_s$ and $\Delta_s$
are homogeneous maps of degree $-s$. By a similar argument as (1),
we may view $m_s \in D^{1,2}_{(-s)}$ and $\Delta_s \in
D^{2,1}_{(-s)}$. Moreover, from (1), we see that $(m_1,\Delta_1)$
can be viewed as an element in $\hat{h}_b^2(B)_{(-1)}$. Now by the
assumption, there exists a homogeneous map $\theta_1:\textsf{m}
\longrightarrow \textsf{m}$ of degree $-1$, such that (see (3.4) and
(3.5))
\begin{align*}
&m_1(a \otimes b)= a \theta_1(b)-\theta_1(ab)+ \theta_1(a)b,\\
&\Delta_1 (c)= \Delta(\theta_1(c))-c_{(1)}\otimes \theta_1(c_{(2)})-
\theta_1(c_{(1)})\otimes c_{(2)}.
\end{align*}
Take $\phi_1 :B[t] \longrightarrow B[t]$ to be a
$\mathbb{K}[t]$-linear map such that
\begin{align*}
\phi_1(a)=a+ \theta_1(a)t,\quad   a \in B.
\end{align*}
Note that $\phi_1$ is a bijective map preserving the identity $1_B$
and the counit $\varepsilon_t$. Consider the deformation
\begin{align*}
(B[t], {m_t}'=\phi_1 \circ m_t \circ (\phi_1^{-1}\otimes
\phi_1^{-1}), {\Delta_t}'=(\phi_1\otimes \phi_1) \circ \Delta_t
\circ \phi_1^{-1}).
\end{align*}
We have
\begin{align*}
{m_t}'(a \otimes b)&=ab + {m_2}'(a \otimes b)t^2+  {m_3}'(a \otimes b)t^3+ \cdots, \\
{\Delta_t}' (c) & = \Delta(c) + {\Delta_2}'(c) t^2 + {\Delta_2}'(c)
t^3+ \cdots
\end{align*}
where ${m_s}'$ and ${\Delta_s}'$ are homogeneous maps of degree
$-s$, $s\geq 2$. Now by comparing (2.3-5) and (3.1-3), we see that
$({m_2}',{\Delta_2}')$ can be viewed as an element in
$\hat{h}_b^2(B)_{(-2)}$. Hence there exists a homogeneous map
$\theta_2:\textsf{m} \longrightarrow \textsf{m}$ of degree $-2$,
such that (again see (3.4) and (3.5))
\begin{align*}
&{m_2}'(a \otimes b)= a \theta_2(b)-\theta_2(ab)+ \theta_2(a)b,\\
&{\Delta_2}' (c)= \Delta(\theta_2(c))-c_{(1)}\otimes
\theta_2(c_{(2)})- \theta_2(c_{(1)})\otimes c_{(2)}.
\end{align*}
Take $\phi_2: B[t] \longrightarrow B[t]$ to be a
$\mathbb{K}[t]$-linear map such that
\begin{align*}
\phi_2(a)=a+ \theta_2(a)t^2,\quad \forall a \in B.
\end{align*}
Now consider the following deformation
\begin{align*}
(B[t], {m_t}''=\phi_2 \circ {m_t}' \circ (\phi_2^{-1}\otimes
\phi_2^{-1}), {\Delta_t}''=(\phi_2\otimes \phi_2) \circ {\Delta_t}'
\circ \phi_2^{-1}),
\end{align*}
whose coefficients of $t$ and $t^2$ vanishes. In other words,
\begin{align*}
{m_t}''(a \otimes b)&=ab + {m_3}''(a \otimes b)t^3+  {m_3}''(a \otimes b)t^4+ \cdots, \\
{\Delta_t}'' (c) & = \Delta(c) + {\Delta_3}''(c) t^3 +
{\Delta_2}''(c) t^4+ \cdots
\end{align*}
Similarly, we may view that  $({m_3}'', {\Delta_3}'')$ lies in
$\hat{h}_b^2(B)_{(-3)}$. By assumption  and comparing (3.1-3), we
have a homogeneous map $\theta_3: \textsf{m}\longrightarrow
\textsf{m}$ such that
\begin{align*}
&{m_3}''(a \otimes b)= a \theta_3(b)-\theta_3(ab)+ \theta_3(a)b,\\
&{\Delta_3}'' (c)= \Delta(\theta_3(c))-c_{(1)}\otimes
\theta_3(c_{(2)})- \theta_3(c_{(1)})\otimes c_{(2)}.
\end{align*}
Now define $\phi_3:B[t] \longrightarrow B[t]$ to be a
$\mathbb{K}[t]$-linear map such that
\begin{align*}
\phi_3(a)=a+ \theta_3(a)t^3,\quad \forall a \in B.
\end{align*}
Thus we get the following deformation
\begin{align*} (B[t],
{m_t}'''=\phi_3 \circ {m_t}'' \circ (\phi_3^{-1}\otimes
\phi_3^{-1}), {\Delta_t}'''=(\phi_3\otimes \phi_3) \circ
{\Delta_t}'' \circ \phi_3^{-1}),
\end{align*}
whose coefficients of $t$, $t^2$ and  $t^3$ vanishes. Now one can
define $\theta_4$ and $\phi_4$, and so on.
\par\vskip 5pt

Finally, define the infinite composition $\cdots \phi_3\circ
\phi_2\circ \phi_1$ to be $\phi$. Note that the
$\mathbb{K}[t]$-linear isomorphism $\phi:B[t] \longrightarrow B[t]$
is well-defined on every $a\in B$, which preserves the identity
$1_B$ and the counit $\varepsilon_t$. (In fact,
$\phi_s(a)=a+\theta_s(a)t^s$ where $\theta_s:\textsf{m}
\longrightarrow \textsf{m}$ is homogeneous of degree $-s$, hence,
for each fixed $a\in B_{(i)}$, $\phi_s(a)=a$ for $s \geq i$.
Consequently, $\phi(a)$ has only nonzero coeffecients  of $t^s$ for
$0 \leq s \leq i$.) By the construction of each map $\phi_s$, we
obtain that the deformation $(B[t], \phi \circ m_t \circ
(\phi^{-1}\otimes \phi^{-1}), (\phi\otimes \phi) \circ \Delta_t
\circ \phi^{-1})$ is trivial, which is also equivalent to the given
deformation. Thus we have proved (2).

\par \vskip 10pt

\noindent (3). \ Let $(B[t]/{(t^{l+1})}, m_t^l, \Delta_t^l)$ be an
element in $\mathcal{E}^l(B)$. Write
\begin{align*}
m_t^l(a \otimes b)=\sum_{0 \leq s \leq l} m_s(a \otimes b) t^s \quad
\mbox{and}\quad  \Delta_t^l(c)=\sum_{0 \leq s \leq l}
\Delta_s(c)t^s,
\end{align*}
where $m_s$ and $\Delta_s$ are homogeneous maps of degree $-s$. By
the same argument as above, one can show that $m_s$ (\emph{resp}.
$\Delta_s$) can be viewed as maps from $\textsf{m}\otimes
\textsf{m}$ to $\textsf{m}$ (\emph{resp}. from $\textsf{m}$ to
$\textsf{m} \otimes \textsf{m}$).\par \vskip 3pt

To extend $(B[t]/{(t^{l+1})}, m_t^l, \Delta_t^l)$ to some element in
$\mathcal{E}^{l+1}(B)$, we just need to find some homogeneous maps
$f:\textsf{m} \otimes \textsf{m} \longrightarrow \textsf{m}$ and $g:
\textsf{m} \longrightarrow \textsf{m}\otimes \textsf{m}$ of degree
$-(l+1)$ such that $(B[t]/{(t^{l+2})}, m_t^l+t^{l+1}f,
\Delta_t^l+t^{l+1}g)$ is an bialgebra over
$\mathbb{K}[t]/{(t^{l+2})}$.\par
 The associativity of $
m_t^l+t^{l+1}f$ is equivalent to
\begin{align*}
  (m_t^l+t^{l+1}f)(((m_t^l+t^{l+1}f)(a \otimes b)) \otimes c))=
  (m_t^l+t^{l+1}f)(a \otimes (( m_t^l+t^{l+1}f)(b\otimes c))),
\end{align*}
for all $a, b , c \in B$. Since $m_t^l$ is associative, then  the
above identity holds if and only if  the two-sides have the same
coefficients of the term $t^{l+1}$. Thus by direct computation, we
get
\begin{align*}
F(a\otimes b \otimes c):&= \sum_{s=1}^l m_s(m_{l+1-s}(a\otimes
b)\otimes c)-m_s(a \otimes m_{l+1-s}(b \otimes c))\\
 &=af(b \otimes c)-f(ab \otimes c)+f(a\otimes bc)-f(a \otimes b)c\\
 &=\delta_{h}^{1,2}(f)(a\otimes b \otimes c).
\end{align*}

Similarly, one obtains that the compatibility of the multiplication
$m_t^l+t^{l+1}f$ and comultiplication $\Delta_t^l+t^{l+1}g$, and the
coassociativity of $\Delta_t^l+t^{l+1}g$ are equivalent to the
following two identities, respectively,

\begin{align*}
H(a\otimes b):&=\sum_{s=1}^l\Delta_s(m_{l+1-s}(a\otimes
b))-\sum_{s+r+s'+r'=l+1}(m_{s'}\otimes m_{r'}) \circ
\tau_{23}(\Delta_s(a)\otimes
\Delta_r(b))\\
&=(\delta_c^{1,2}(f)+\delta_h^{2,1}(g))(a \otimes b),
\end{align*}
and
\begin{align*}
G(c):&= \sum_{s=1}^l(\Delta_s \otimes {\rm Id})\circ
\Delta_{l+1-s}(c)-({\rm Id}\otimes \Delta_s)\circ\Delta_{l+1-s}(c)\\
&=c_{(1)}\otimes g(c_{(2)})-(\Delta\otimes {\rm Id})(g(c))+({\rm
Id}\otimes \Delta)(g(c))-g(c_{(1)})\otimes c_{(2)}\\
&=\delta_c^{2,1}(g)(c),
\end{align*}
where $a, b, c \in \textsf{m}$, and $\tau_{23}$ is the flip map with
respect to the second and third positions.\par

Now it is direct to check that the element $(F,H,G)\in
\hat{D}_{(-l-1)}$ is a cocycle (exactly as  \cite{G} in the case of
algebras and \cite{GS2} in the case of non-graded bialgebras), i.e.,
it lies in the kernel of the differential $\partial_{(-l-1)}^3$ from
(2.3-5), therefore, it can be viewed as an element in the cohomology
group $\hat{h}_b^3(B)_{(-l-1)}$. Now by comparing the above three
identities with (3.6), we obtain that if
$\hat{h}_b^3(B)_{(-l-1)}=0$, then such maps $f, g$ always exist,
i.e., we can extend $(B[t]/{(t^{l+1})}, m_t^l, \Delta_t^l)$ to
$(B[t]/{(t^{l+2})}, m_t^l+t^{l+1}f, \Delta_t^l+t^{l+1}g)$. Note that
by the above three equivalences, one deduces that
$(B[t]/{(t^{l+2})}, m_t^l+t^{l+1}f, \Delta_t^l+t^{l+1}g)$ belongs to
$\mathcal{E}^{l+1}(B)$. This completes the proof.\hfill
$\blacksquare$

\vskip 10pt

\noindent{\bf Acknowledgement:}\quad We thank the referees for many
valuable suggestions.

\bibliography{}

\end{document}